\documentclass[aps,prl,twocolumn,superscriptaddress]{revtex4-1}

\usepackage{epsfig}\usepackage{graphicx}
\usepackage{amsmath,amssymb}
\usepackage{amsthm}
\newcommand{\comment}[1]{}
\usepackage{xcolor}
\newcommand{\peq}{p_{{\rm eq}}}
\newtheorem{cor}{Corollary}[]

\begin{document}


\title{Nonparametric forecasting of low-dimensional dynamical systems}


\author{Tyrus Berry}
\affiliation{Department of Mathematics, the Pennsylvania State University, University Park, PA 16802-6400, USA}
\author{Dimitrios Giannakis}
\affiliation{Courant Institute of Mathematical Sciences, New York University, New York, NY 10012, USA}
\author{John Harlim}
\affiliation{Department of Mathematics, the Pennsylvania State University, University Park, PA 16802-6400, USA}
\affiliation{Department of Meteorology, the Pennsylvania State University, University Park, PA 16802-5013, USA}


\date{\today}

\begin{abstract}
This letter presents a non-parametric modeling approach for forecasting stochastic dynamical systems on low-dimensional manifolds. The key idea is to represent the discrete shift maps on a smooth basis which can be obtained by the diffusion maps algorithm. In the limit of large data, this approach converges to a Galerkin projection of the semigroup solution to the underlying dynamics on a basis adapted to the invariant measure. This approach allows one to quantify uncertainties (in fact, evolve the probability distribution) for non-trivial dynamical systems with equation-free modeling. We verify our approach on various examples, ranging from an inhomogeneous anisotropic stochastic differential equation on a torus, the chaotic Lorenz three-dimensional model, and the Ni\~{n}o-3.4 data set which is used as a proxy of the El-Ni\~{n}o Southern Oscillation.
\end{abstract}

\maketitle


A significant challenge in modeling is to account for physical processes which are often not well understood but for which large data sets are available. A standard approach is to perform regression fitting of the data into various parametric models. While this approach is popular and successful in many applied domains, the resulting predictive skill can be sensitive to the choice of models, the parameter fitting algorithm, and the complexity of the underlying physical processes. An alternative approach is to avoid choosing particular models and/or parameter estimation algorithms and instead apply a nonparametric modeling technique. In particular, nonparametric modeling based on local linearization of discrete shift maps on paths has been successful in predicting mean statistics of data sets generated by dynamical systems with low-dimensional attractors \cite{shiftop1,shiftop2,shiftop3}. 

In this letter, we generalize this nonparametric approach to quantify the evolving probability distribution of the underlying dynamical system. The key idea is to project the shift map on a set of basis functions that are generated by the diffusion maps algorithm \cite{diffusion} with a variable bandwidth diffusion kernel \cite{BH14VB}.  This approach has connections with a recently developed family of kernels \cite{G14Cone}, which utilize small shifts of a deterministic time series to estimate the dynamical vector field. The method also generalizes a recently introduced non-parametric modeling framework for gradient systems \cite{BH14UQ} to inhomogeneous stochastic systems having non-gradient drift and anisotropic diffusion.  Consider a dynamical system, 
\begin{align}\label{SDE} dx = a(x)\, dt + b(x)\, dW_t \end{align}
where $W_t$ is a standard Brownian process, $a(x)$ a vector field, and $b(x)$ a diffusion tensor, all defined on a manifold $\mathcal{M}\subset \mathbb{R}^n$. Given a time series $x_i = x(t_i)$, generated by \eqref{SDE} at discrete times $\{t_i\}_{i=1}^{N+1}$, we are interested in constructing a forecast model so that given an initial density $p_0(x)$ we can estimate the density $p(x,t)= e^{t\mathcal{L}^*}p_0(x)$ at time $t>0$, without the Fokker-Planck operator, $\mathcal{L}^*$, of \eqref{SDE} and without knowing or estimating $a$ and $b$. We will assume that \eqref{SDE} is ergodic so that $\{x_i\}$ are sampled from the invariant measure $\peq(x)$ of \eqref{SDE}.  Note that all probability densities are relative to a volume form $dV$ which $\mathcal{M}$ inherits from the ambient space, and the generator $\mathcal{L}$ of \eqref{SDE} is the adjoint of $\mathcal{L}^*$ with respect to $L^2(\mathcal{M},dV)$.

We note that our approach differs significantly from previous approaches such as \cite{PRL1,PRL2}, which estimate $a$ and $b$ explicitly in the ambient space $\mathbb{R}^n$, relying on the Kramer-Moyal expansion.  In contrast, our approach directly estimates the semi-group solution $e^{\tau\mathcal{L}}$ on the manifold $\mathcal{M}\subset \mathbb{R}^n$, so that $a$ and $b$ are represented implicitly. The advantages of our approach are that the data requirements are independent of the ambient space dimension and only depend on the intrinsic dimension of $\mathcal{M}$, and we will be able to estimate the semi-group solution, $e^{\tau\mathcal{L}}$, directly from the data for any sampling time $\tau$.  

Our approach is motivated by a rigorous connection between the shift map, $S$ which we define by $Sf(x_i) = f(x_{i+1})$ for a function $f\in L^2(\mathcal{M},\peq)$, and the semi-group solution, $e^{\tau\mathcal{L}}$, of the underlying dynamical system \eqref{SDE}.  Applying It\^o's Lemma to $f(x)$, one can show,
\begin{align}\label{shiftexpansion} 
Sf(x_i) &= e^{\tau\mathcal{L}}f(x_i) + \int_{t_i}^{t_{i+1}}\nabla f^\top b  \, dW_s +  \int_{t_i}^{t_{i+1}} B f\,ds ,
\end{align}
where $\tau = t_{i+1}-t_i$, $Bf = \mathcal{L}f - \mathbb{E}\big[\mathcal{L}f]$ and the expectation is with respect to paths of \eqref{SDE} conditional to $x(t_i)$. The detailed derivation of \eqref{shiftexpansion} is in Appendix B in the supplementary material. Since $\mathbb{E}[Sf(x_i)] = e^{\tau\mathcal{L}}f(x_i)$, we can use the shift map, $S$, to directly estimate the semigroup solution, $e^{\tau\mathcal{L}}$. However, $S$ is a noisy estimate of $e^{\tau\mathcal{L}}$ and our key contribution is to minimize the error  by representing $S$ on a basis of smooth functions. 

Minimizing the error requires choosing a basis which minimizes the functional $\|\nabla f\|_{\peq}$, which is shown in Appendix B. 
Intuitively, this is because we want to bound the stochastic integral in \eqref{shiftexpansion}, whose integrand contains $\nabla f$.
The functional $\|\nabla f\|_{\peq}$ is minimized by the eigenfunctions of the generator, $\mathcal{\hat{L}}$, of a stochastically forced gradient flow with potential $U(x) = -\log(\peq (x))$. Let $\lambda_j$ and $\varphi_j$ be the eigenvalues and eigenfunctions of $\mathcal{\hat{L}}$; $\{\varphi_j\}$ are orthonormal on $L^2(\mathcal{M},\peq)$. Since $\hat{\mathcal{L}}$ is the generator of a gradient flow systems, it is easy to check that $\psi_j=\peq\varphi_j$ are the eigenfunctions of the adjoint operator $\mathcal{\hat{L}}^*$, which are orthonormal on $L^2(\mathcal{M},\peq^{-1})$. Numerically, we obtain $\varphi_j(x_i)$ as eigenvectors of a stochastic matrix, constructed by evaluating a variable bandwidth kernel on all pairs of data points and then applying an appropriate normalization \cite{BH14VB}. We summarize this procedure and the related results in Appendix A in the Supplementary material.  We emphasize that the operator $\hat{\mathcal{L}}$ is used only to estimate $\{\varphi_j\}$ which is the optimal basis for smoothing the shift operator $S$ approximating semi-group solution $e^{\tau\mathcal{L}}$ of the full system \eqref{SDE}.  

We write the solution, $p(x,\tau) = e^{\tau\mathcal{L}^*}p_0(x)$, as follows,
\begin{align}\nonumber
p(x,\tau) = \sum_{l=1}^{\infty} \langle e^{\tau\mathcal{L}^*}p_0,\psi_l \rangle_{\peq^{-1}} \psi_l(x) = \sum_{l=1}^{\infty} \langle p_0,e^{\tau\mathcal{L}}\varphi_l \rangle \psi_l(x).
\end{align}
Similarly, the initial density is $p_0(x) = \sum_j c_j(0) \psi_j(x)$, where $c_j(0) = \langle p_0,\psi_j\rangle_{\peq^{-1}}$. We therefore obtain,
\begin{align} \label{p2}
p(x,\tau) = \sum_{l=1}^{\infty}\sum_{j=1}^{\infty} c_j(0) A_{lj}(\tau)\, \peq(x)\varphi_l(x),
\end{align}
where $A_{lj}(\tau):=\langle\varphi_j,e^{\tau\mathcal{L}}\varphi_l \rangle_{\peq}$. Based on the discussion after \eqref{shiftexpansion}, we will use $\langle\varphi_j,S\varphi_l \rangle_{\peq}$ to estimate $A_{lj}$. Numerically, we estimate $A_{lj}$ by a Monte-Carlo integral,
\begin{align}
\hat{A}_{lj} &= \frac{1}{N}\sum_{i=1}^{N} \varphi_j(x_i)\varphi_l(x_{i+1}),\label{AjlMC}
\end{align}
such that, the {\it diffusion forecast} is defined as follows:
\begin{align} p(x,\tau) \approx \sum_{l=1}^{\infty} \peq(x) \,\varphi_l(x) \sum_{j=1}^{\infty} \hat{A}_{lj}(\tau)c_j(0).\label{papprox}\end{align}

One can show that $\mathbb{E}[\hat{A}_{lj}]=A_{lj}$ which means that $\hat{A}_{lj}$ is an unbiased estimate of $A_{lj}$. Moreover, the error of this estimate is of order $\lambda_l\sqrt{\tau/N}$ in probability assuming that $x_i$ are independent samples of $\peq$. This shows that we can apply a diffusion forecast for any sampling time $\tau$ given a sufficiently large data set $N$.  For more details, see Appendix B in the Supplementary Material.

{\it Non-gradient drift anisotrophic diffusion:} We first verify the above approach for a system of SDE's of the form \eqref{SDE} on a torus defined in the intrinsic coordinates $(\theta,\phi) \in [0,2\pi)^2$ with drift and diffusion coefficients,
\begin{align}
a(\theta,\phi) &= \begin{pmatrix}\frac{1}{2} + \frac{1}{8}\cos(\theta)\cos(2\phi)+\frac{1}{2}\cos(\theta+\pi/2) \\ 10+\frac{1}{2}\cos(\theta+\phi/2)+\cos(\theta+\pi/2))\end{pmatrix},\nonumber\\
b(\theta,\phi) &= \begin{pmatrix} \frac{1}{4}+\frac{1}{4}\sin(\theta) & \frac{1}{4}\cos(\theta + \phi) \\ \frac{1}{4}\cos(\theta+\phi) & \frac{1}{40}+ \frac{1}{40}\sin(\phi)\cos(\theta) \end{pmatrix}. \nonumber
\end{align}
This example is chosen to exhibit non-gradient drift, anisotropic diffusion, and multiple time scales.  Since it is a system of the form \eqref{SDE} on a smooth manifold, our theory shows that the shift operator $S$ is an unbiased estimator for the semigroup solution $e^{\tau\mathcal{L}}$, and in the limit of large data the diffusion forecast will capture all aspects of the evolution of the density $p(x,t)$.  

We now verify this theory using a training data set of 20000 points generated by numerically solving the SDE in \eqref{SDE} with a discrete time step $\Delta t = 0.1$ and then mapping this data into the ambient space, $\mathbb{R}^3$, via the standard embedding of the torus given by $(x,y,z) = ((2+\sin(\theta))\cos(\phi),(2+\sin(\theta))\sin(\phi),\cos(\theta))$.  We define a Gaussian initial density $p_0(\theta,\phi)$ with a randomly selected mean and a diagonal covariance matrix with variance 1/10. The initial density is projected into a basis of $M=1000$ eigenfunctions giving coefficients $c_j(0) = \langle p_0/\peq,\varphi_j  \rangle_{\peq} \approx \sum_{i=1}^{20000} p_0(\theta_i,\phi_i)\varphi_j(\theta_i,\phi_i)/\peq(\theta_i,\phi_i)$.  The coefficients $c_j(0)$ are evolved forward in time in discrete steps of length $\Delta t = 0.1$ by $\hat{A}$, constructed by \eqref{AjlMC} so that the forecast at time $\tau=n\Delta t$ is effectively $\hat{A}(\Delta t)^n$.

\begin{figure}[htbp]
\centering
\includegraphics[height=.35\textwidth,width=0.48\textwidth]{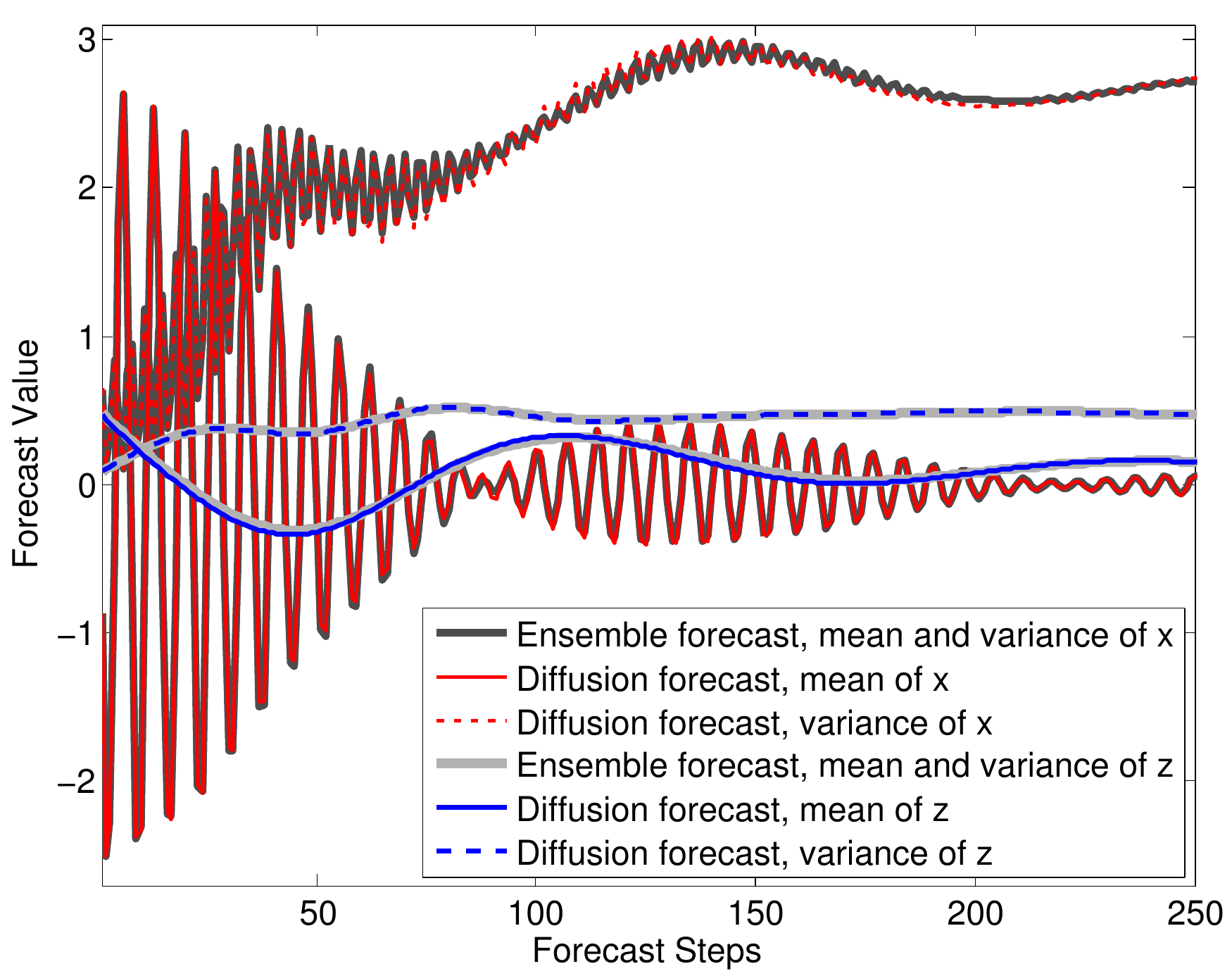}
\caption{\label{TorusEx} Validation of first two moments of the diffusion forecast for a stochastic dynamical system on a torus in $\mathbb{R}^3$.}\label{example1}
\end{figure}

In Figure~\ref{example1}, we show the evolution of the first two-moments in the ambient space, for the fast and slow variables, $x$ and $z$, respectively, created by the diffusion forecast in \eqref{papprox}. To verify the accuracy of the diffusion forecast, we also show the corresponding moments produced by an ensemble forecast of 50000 initial conditions, randomly sampled from the initial distribution $p_0(\theta,\phi)$, evolved using the true dynamical system. Notice the long-time pathwise agreement of both moments constructed via the diffusion forecast and those constructed by an ensemble forecast. See also a video of the evolution of $p$ in the Supplementary Material.  

{\it Lorenz-63 model:}
Next, we apply the diffusion forecasting algorithm to data generated by the Lorenz-63 model \cite{L63}, with error in the initial state.  Unlike the previous example, this model does not technically satisfy the requirements of our theory because the attractor is a fractal set rather than a smooth manifold.  Our results indicate that the applicability of the method seems to extend beyond the current theory.  We compare our approach to the classical nonparametric prediction method which uses a local linear approximation of the shift map \cite{shiftop1,shiftop2,shiftop3} and a standard ensemble forecasting method with the true model, applied with 50000 initial conditions, sampled from the same initial distribution $p_0$.  

We generate 10000 data points with discrete time steps $\Delta t = 0.1$ and $\Delta t = 0.5$ by solving the Lorenz-63 model. We use the first 5000 data points as training data for the nonparametric models, and the remaining 5000 data points to verify the forecasting skill. For each of the 5000 verification points, $x_t$, we define the initial state $\hat x_t = x_t + \xi_t$ by introducing random perturbations $\xi_t$ sampled from $\mathcal{N}(0,0.01)$.  Each forecast method starts with the same initial density, $p_0 = \mathcal{N}(\hat x_t,0.01)$ centered at the perturbed verification point. We chose this very small perturbation to demonstrate the diffusion forecast for an initial condition which is almost perfect; as the amount of noise increases the advantage of the diffusion forecast over the linear methods is even more significant.  

\begin{figure}[htbp]
\centering
\includegraphics[height=.285\textwidth,width=0.48\textwidth]{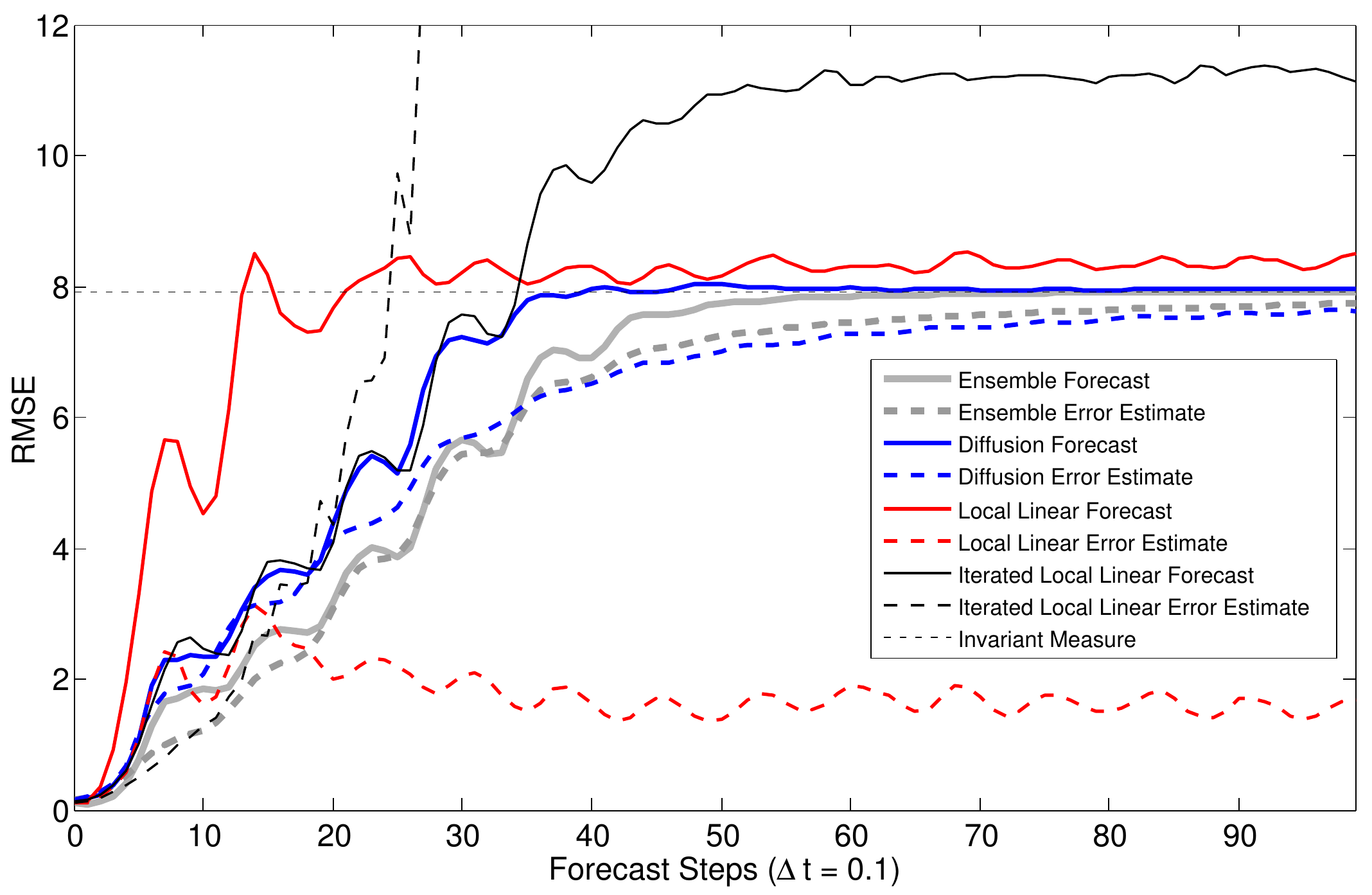}
\includegraphics[height=0.285\textwidth,width=0.48\textwidth]{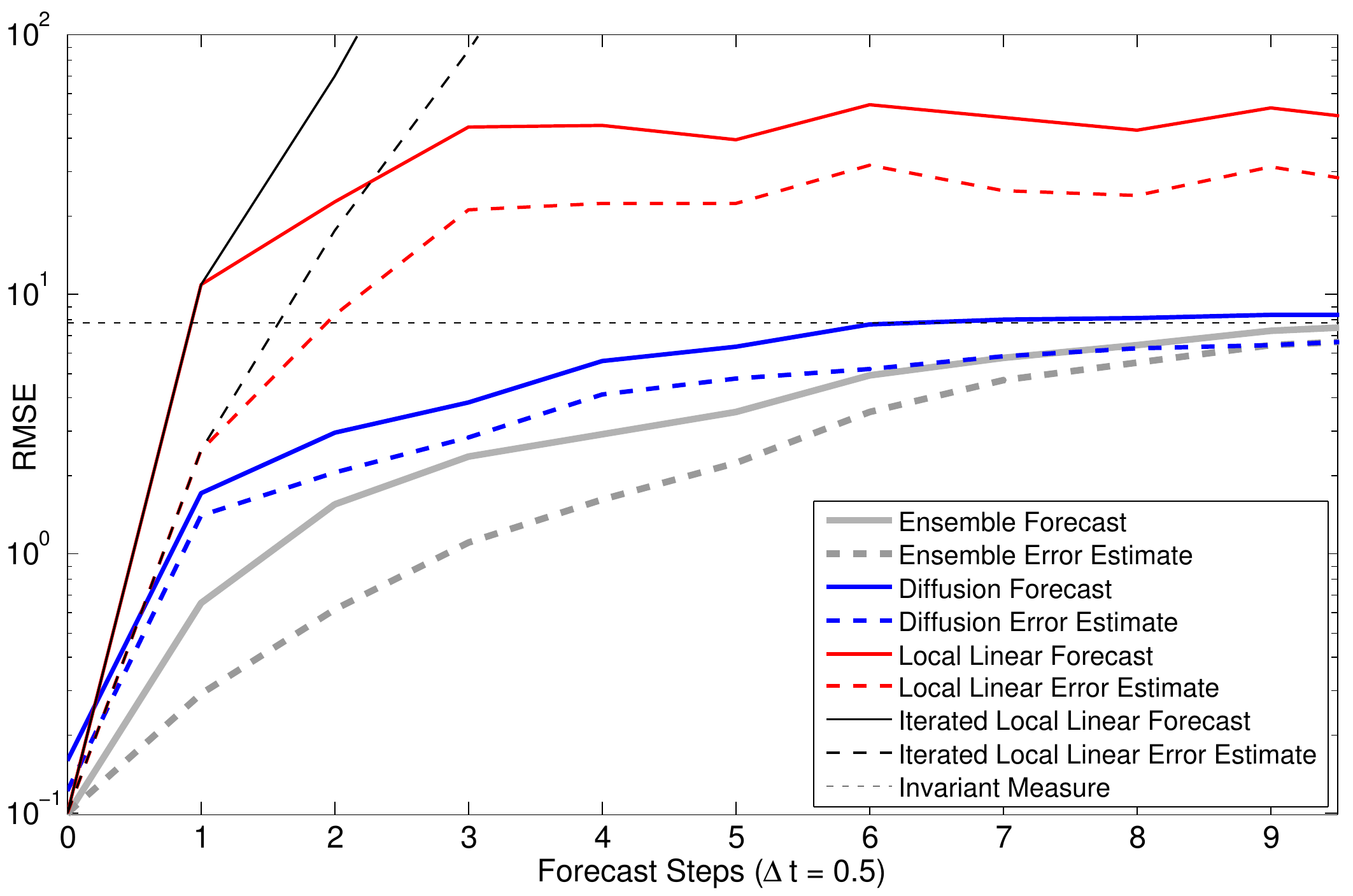}
\includegraphics[height=0.38\textwidth,width=0.48\textwidth]{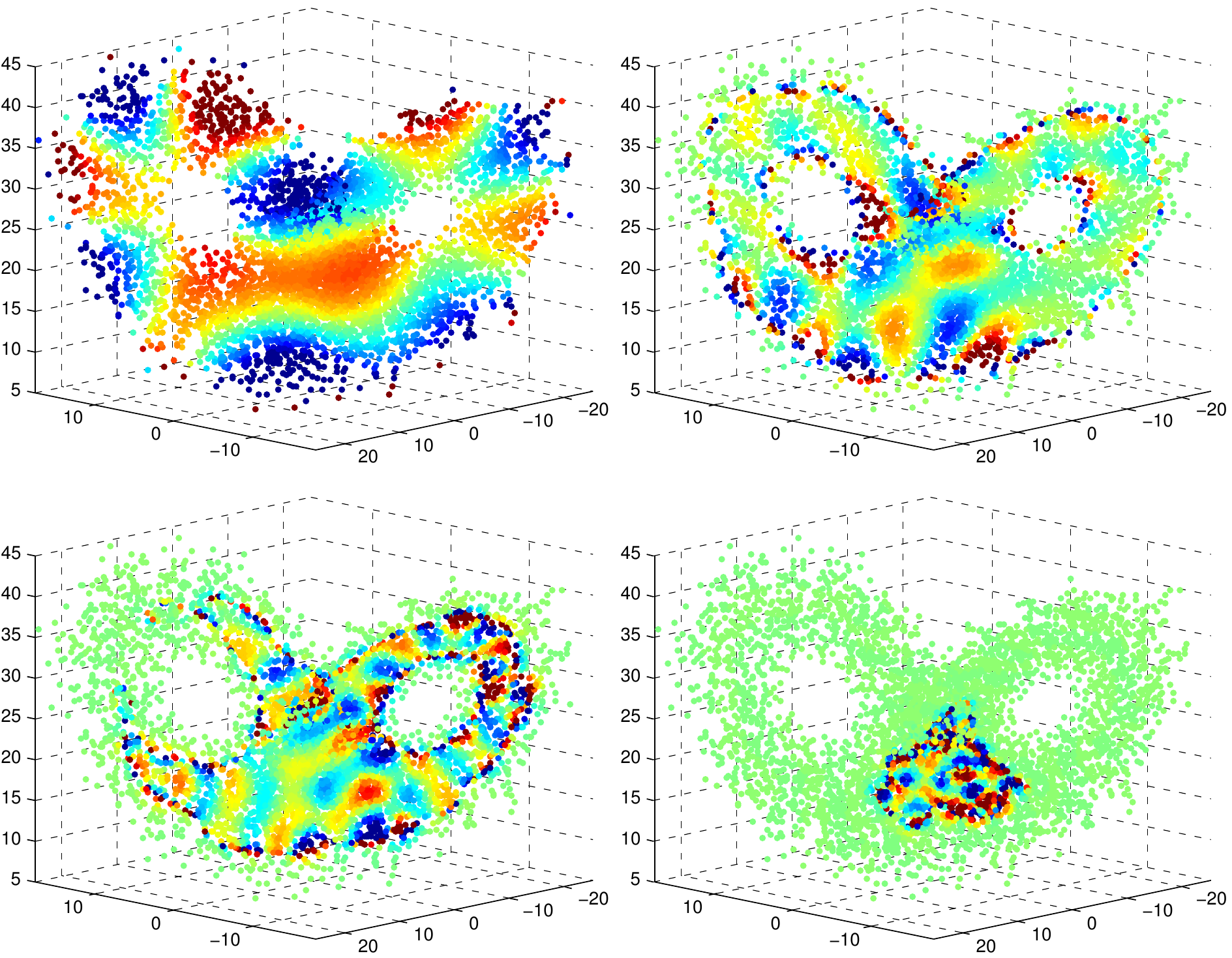}
\caption{\label{L63Ex} Comparing forecast methods for Lorenz-63, Top: $\Delta t = 0.1$, Middle: $\Delta t = 0.5$. Bottom: Eigenfunctions $\varphi_{40},\varphi_{500},\varphi_{1500},$ and $\varphi_{4000}$ of the coarse approximation of the fractal attractor by a manifold.}
\vspace{-15pt}
\end{figure}

The diffusion forecast is performed with 4500 eigenfunctions $\varphi_j$, constructed with the diffusion maps algorithm with a variable bandwidth \cite{BH14VB} (we show examples in Figure \ref{L63Ex}). The local linear forecast uses ordinary least squares to fit an affine model to the $n$-step shift map on the 15 nearest neighbors to the initial state. The iterated local linear forecast completes this process for one step and then recomputes the 15 nearest neighbors to the 1-step forecast and then repeats the process. The variance estimate of the local linear models is given by conjugating the covariance matrix of $p_0$ with the linear part of the appropriate affine forecast model. We compute the root mean squared error (RMSE) between each mean forecast and the true state, averaged over the verification period of 5000 steps.  We also show the standard deviation of the forecast density, so that a forecasting method has good uncertainty quantification (UQ) if the standard deviation agrees with the RMSE.  

Of course, the ensemble forecast with the true model gives the best forecast, however the diffusion forecast is a considerable improvement over the local linear forecast.  For short $\Delta t = 0.1$, the iterated local linear forecast is comparable to the diffusion forecast except in the long term where the iterated local linear forecast exhibits significant bias.  Moreover, the iterated local linear forecast significantly overestimates the error variance in the intermediate to long term forecast.  This overestimation is due to the positive Lyapunov exponent, which is implicitly estimated by the product of the iterated local linearizations.  In contrast, the direct local linearization is unbiased in the long term, but converges very quickly to the invariant measure and underestimates the variance.  This underestimation is because no single linearization can capture the information creation introduced by the positive Lyapunov exponent.  For long $\Delta t = 0.5$, the bias in the local linear models leads them to diverge far beyond the invariant measure for even intermediate term forecasts. The ensemble forecast provides the most consistent UQ since it has access to the true model, however the diffusion forecast produces reasonable estimates without knowing the true model as shown in Figure \ref{L63Ex}.  The local linear forecast error estimates vary widely and do not robustly provide a useful UQ whereas the diffusion forecast is robust across multiple sampling times as suggested by the theory. We include a video showing good long term agreement between the diffusion forecast density and an ensemble in the Supplementary Material.

The diffusion forecast is able to give a reasonable estimate of the evolution of the density by building a consistent finite dimensional Markovian approximation of the dynamics.  This Markovian system incorporates global knowledge of the attractor structure via the smoothing with the adapted basis $\{\varphi_j\}$.  This Markovian approximation of the Lorenz-63 model implicitly uses a small Brownian forcing to replicate the entropy generation of the positive Lyapunov exponent.  

\begin{figure}
\centering
\includegraphics[height=0.28\textwidth,width=0.48\textwidth]{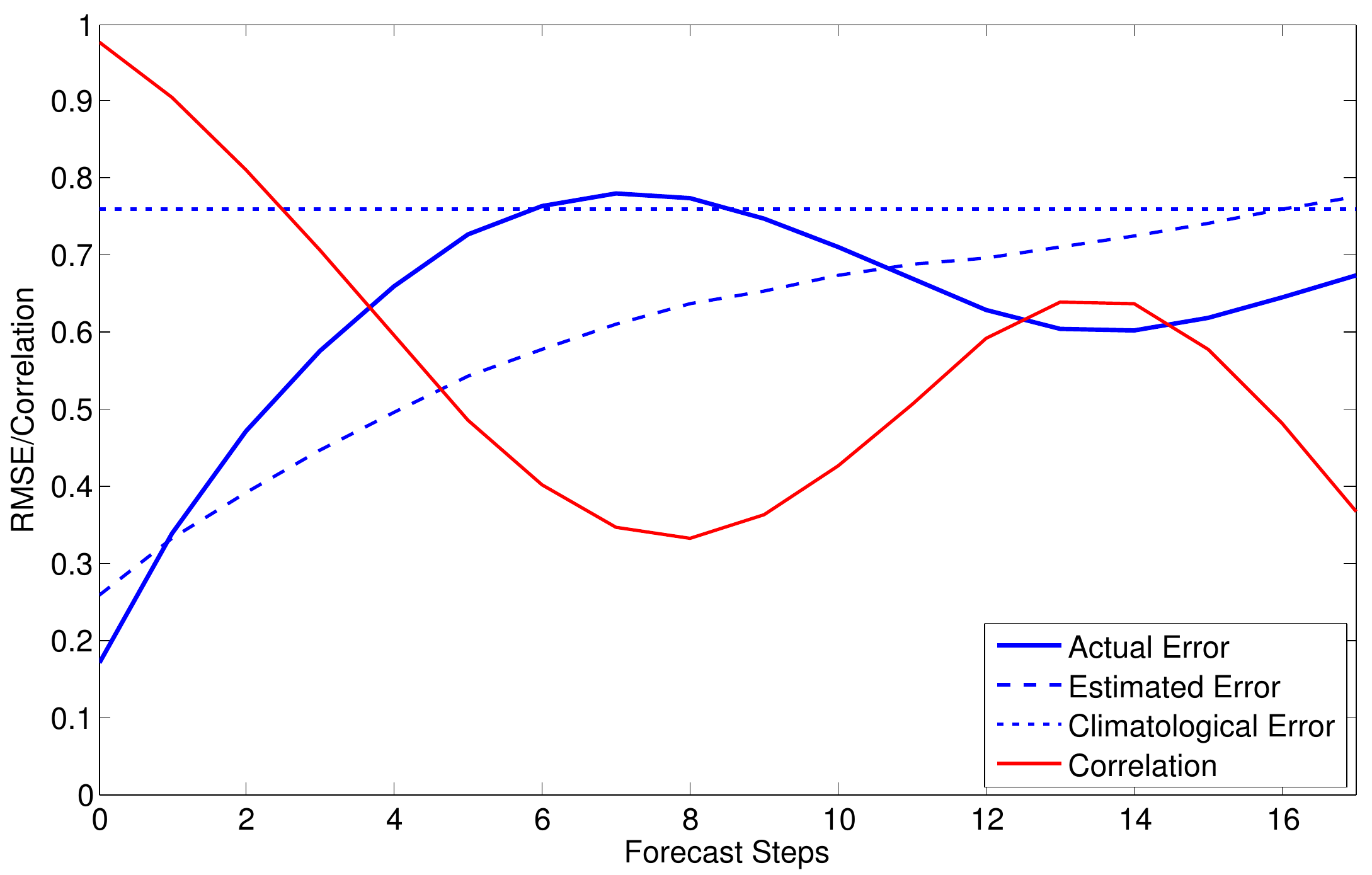}
\includegraphics[height=0.28\textwidth,width=0.48\textwidth]{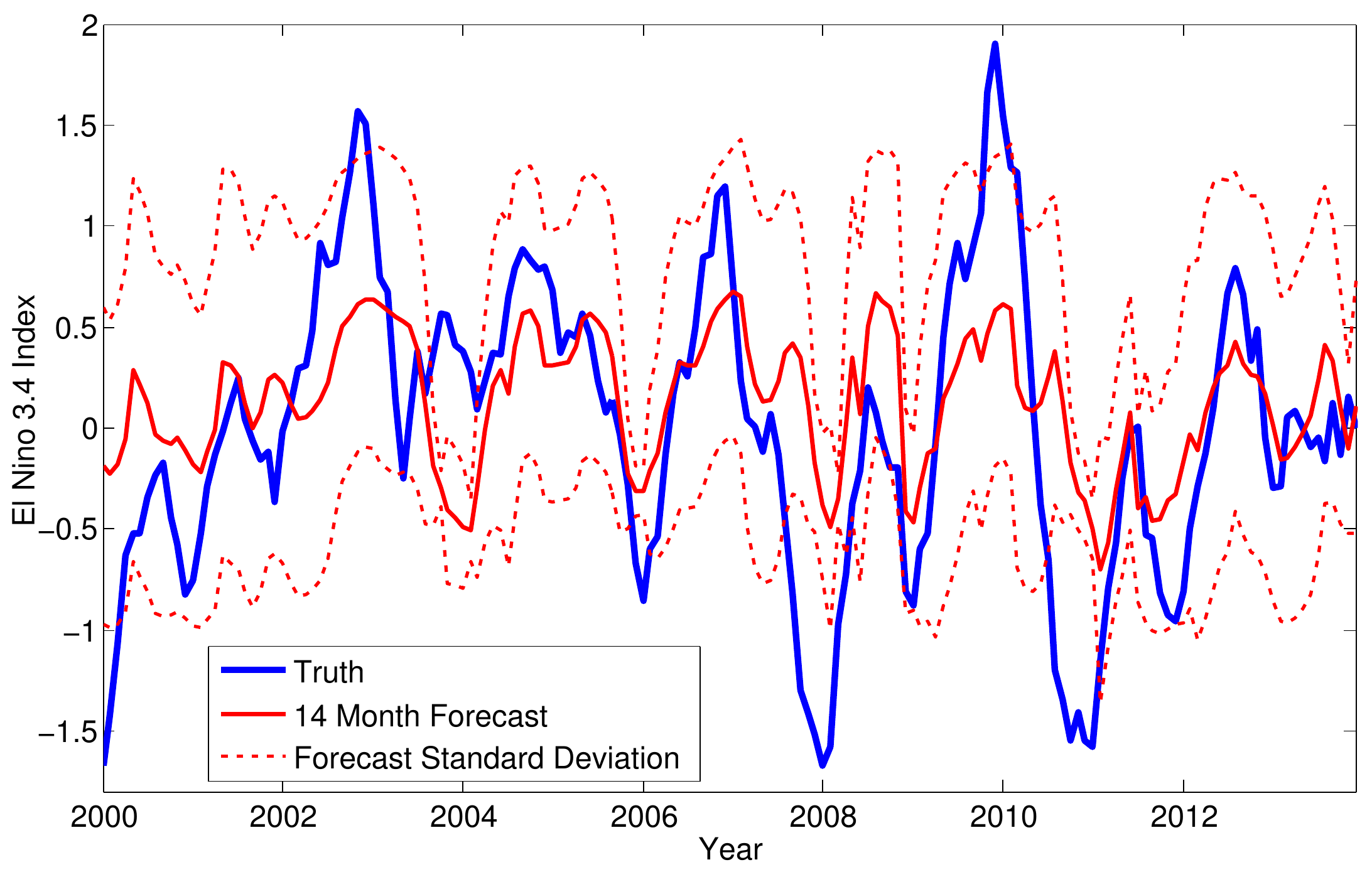}\vspace{-5pt}
\caption{\label{NinoEx} Forecasting for the El Ni\~{n}o 3.4 index. RMS and correlation (top); 14-month lead forecast (bottom)}\label{fig3}\vspace{-15pt}
\end{figure}

{\it El-Ni\~{n}o data set:} We now apply our method to a real world data set, where the validity of our theory is unverifiable, namely the Ni\~{n}o-3.4 index, which records the monthly anomalies of sea surface temperature (SST) in the central equatorial Pacific region (the raw dataset is available from NOAA \url{ftp://ftp.cpc.ncep.noaa.gov/wd52dg/data/indices/}). In applying this method, we implicitly assume that there is an underlying dynamics on a low-dimensional manifold that generates these SST anomalies. Since the time series is one-dimensional, we apply the time-delay embedding technique to the data set, following \cite{SYC,GiannakisPNAS,DMDC}, which will recover this low dimensional manifold if it exists. We use a 5-lag embedding, empirically chosen to maximize the forecast correlation skill between the true time series and the mean estimate.  We construct the $A_{lj}$ matrix with 80 eigenfunctions obtained by the diffusion maps algorithm with a variable bandwidth kernel, applied on the lag-embedded data set.  

In this experiment, we train our nonparametric model with monthly data between Jan 1950-Dec 1999 (600 data points), following \cite{ckg:11} and we verify the forecasting skill on the period of Jan 2000-Sept 2013. The initial distribution $p_0$ is generated with the same method as in the Lorenz-63 example. Based on the RMSE and correlation measure (see Figure~\ref{fig3}), the forecasting skill decays to the climatological error in about 6 months but then the skill improves, peaking at 13-14 month lead time. In fact, our 14-month lead forecast skill, in terms of RMSE 0.60 and correlation 0.64, is significantly better than that of the method proposed in \cite{ckg:11} (Fig.~3 in their paper suggests RMSE 1.4 and correlation 0.4) who claimed to beat the current operational forecasting skill. The 14-month lead forecast mean estimate gives a reasonably correct pattern, and the diffusion forecast provides a reasonable error bar (showing one stdev.) which is a useful UQ. We include a movie in the Supplementary Material showing the nontrivial evolution of the forecast distribution starting 14 months before January 2004.  

Difficulty in improving the forecasts in this problem may be due to combinations of many factors, including the validity of our assumption of the existence of low-dimensional structures, a transient in the time series, a large stochastic component, and memory effects. One possibility is to combine the local linear models via the diffusion basis to form a global model which respects to invariant measure.  Another issue is that both the observational and dynamical noise in the data is currently treated as part of the process, and it would be advantageous to isolate the attractor and build the basis there. Finally, the empirical success of this method suggests that it is possible to approximate a fractal attractor with a smooth manifold, however, there is limited theoretical interpretation for such an approximation. 

This research is supported under the ONR MURI grant N00014-12-1-0912. D.G.\ also acknowledges financial support from ONR DRI grant N00014-14-1-0150.  J.H. also acknowledges financial support from ONR grant N00014-13-1-0797 and the NSF grant DMS-1317919. 

\bibliography{VBbib}

\section{Appendix A: Estimating Elliptic Operators with Diffusion Kernels}

In this appendix we include the details of the numerical algorithm introduced in \cite{BH14VB} which is used to estimate the eigenfunctions $\varphi_j$ of the elliptic operator $\hat{\mathcal{L}} = \Delta - c_1 \nabla U \cdot \nabla$. We assume that the data set $\{x_i\}_{i=1}^N \subset \mathbb{R}^n$ lies on a $d$-dimensional submanifold $\mathcal{M}$ of the ambient Euclidean space $\mathbb{R}^n$.  In order to describe functions on the manifold, we will represent a function $f$ by evaluation on the data set, so that $f$ is represented by the discrete $N\times 1$ vector $\vec f = (f_1,...,f_N)^\top$ where $f_i = f(x_i)$.  In this framework, operators which map functions to functions are represented as $N\times N$ matrices, since these take the $N\times 1$ functions to $N\times 1$ functions.  In this sense, we will approximate the eigenfunctions of $\hat{\mathcal{L}}$ as eigenvectors of a kernel matrix which provably approximates the operator $\hat{\mathcal{L}}$ in the limit of large data.  In this paper we are interested in data sets generated by dynamical systems, however the general theory of this section is valid for arbitrary data sets.  In this case, we will assume that the dynamics are ergodic so that the manifold $\mathcal{M}$ is an attractor of the system, and the sampling measure $q$ of the data on $\mathcal{M}$ is the same as the invariant measure $\peq$ of the dynamics.

We first note that $\hat{\mathcal{L}}$ is the generator of the gradient flow system, 
\begin{align}\label{gradSDE} dx = -c_1 \nabla U(x) \, dt + \sqrt{2} \, dW_t, \end{align}
where $W_t$ is a Brownian motion on the manifold $\mathcal{M}$ (so that the Laplacian $\Delta$ is the infinitesimal generator of $W_t$).  The potential function $U = -\log(q)$ is determined by the sampling measure, $q$.  In particular, we will be interested in the case $c_1 =1$, in which case the invariant measure of the system \eqref{gradSDE} is $e^{-c_1 U} = q = \peq$.  This means that the invariant measure of the true dynamical system which governs the evolution of the data set is the same as that of the gradient flow system generated by $\hat{\mathcal{L}}$.  Since $\hat{\mathcal{L}}$ is a negative semi-definite elliptic operator, self-adjoint with respect to $\peq$, the eigenfunctions $\{\varphi_j\}$ form a basis for $L^2(\mathcal{M},\peq)$.  Moreover, this is the smoothest basis with respect to $\peq$ in the sense that each $\varphi_j$ minimizes the norm $||\nabla f||_{\peq}= -\langle f,\hat{\mathcal{L}}f \rangle_{\peq}$ subject to $\varphi_j$ being orthogonal to all $\{\varphi_l\}_{l<j}$.  The minimal value of the norm is given by the corresponding eigenvalue, $||\nabla \varphi_j||_{\peq} = \lambda_j$.

Approximating the operator $\hat{\mathcal{L}}$ was first achieved in \cite{diffusion} for compact manifolds using a fixed bandwidth kernel.  However, in order to allow the manifold $\mathcal{M}$ to be non-compact, it was shown in \cite{BH14VB} that we must use a variable bandwidth diffusion kernel of the form,
\begin{align}\label{vbkernel} K^S_{\epsilon}(x,y) = \exp\left(-\frac{||x-y||^2}{4\epsilon (q_{\epsilon}(x)q_{\epsilon}(y))^{\, \beta}}\right), \end{align} 
where $q_{\epsilon}(x) = q(x) + \mathcal{O}(\epsilon)$ is an order-$\epsilon$ estimate of the sampling density and $\beta$ should be negative so that the bandwidth function $q_{\epsilon}(x)^{\beta}$ is large in regions of sparse sampling and small in regions of dense sampling.  The algorithm presented below is closely related to that presented in \cite{BH14VB} and is motivated by Corollary 1 of \cite{BH14VB}.

\begin{cor}\label{corollary}  Let $q\in L^1(\mathcal{M}) \cap \mathcal{C}^3(\mathcal{M})$ be a density that is bounded above on an embedded $d$-dimensional manifold $\mathcal{M} \subset \mathbb{R}^n$ without boundary and let $\{x_i\}_{i=1}^N$ be sampled independently with distribution $q$.  
Let $K^S_{\epsilon}$ be a variable bandwidth kernel of the form \eqref{vbkernel} with bandwidth function $q_{\epsilon}^{\beta}$ where $q_{\epsilon} = q + \mathcal{O}(\epsilon)$ is any order-$\epsilon$ estimate of $q$. For a function $f \in L^2(\mathcal{M},q)\cap \mathcal{C}^3(\mathcal{M})$ and an arbitrary point $x_i \in \mathcal{M}$, define the discrete functionals,
\begin{align}
F_i(x_j) &= \frac{K^S_{\epsilon}(x_i,x_j)f(x_j)}{q^S_{\epsilon}(x_i)^\alpha q^S_{\epsilon}(x_j)^\alpha},\quad
G_i(x_j) = \frac{K^S_{\epsilon}(x_i,x_j)}{q^S_{\epsilon}(x_i)^\alpha q^S_{\epsilon}(x_j)^\alpha}, \nonumber
\end{align}
where $q^S_{\epsilon}(x_i) = \sum_l K^S_{\epsilon}(x_i,x_l)/q_{\epsilon}(x_i)^{d\beta}$.  Then,
\begin{align}\label{errorestimate}
&L^S_{\epsilon,\alpha,\beta}f(x_i)  \equiv  \frac{1}{\epsilon m q_{\epsilon}(x_i)^{2\beta}}\left(\frac{\sum_{j}F_i(x_j)}{\sum_{j}G_i(x_j)}-f(x_i)\right)  \\
&\hspace{2pt}= \hat{\mathcal{L}} f(x_i) + \mathcal{O}\left(\epsilon, \frac{q(x_i)^{(1-d\beta)/2}}{\sqrt{N}\epsilon^{2+d/4}} ,\frac{||\nabla f(x_i)||q(x_i)^{-c_2}}{\sqrt{N}\epsilon^{1/2+d/4}} \right), \nonumber
\end{align}
with high probability, where $c_1 = 2-2\alpha + d\beta + 2\beta$ and $c_2 = 1/2-2\alpha+2d\alpha +d\beta/2+\beta$ and $m$ is a constant depending on the form of the kernel, and $m=2$ for \eqref{vbkernel}.
\end{cor}

The key result of \cite{BH14VB} is that for the error to be bounded when $q$ is not bounded below, we require $c_2<0$ (otherwise as $q\to 0$ the final error term becomes unbounded).  Since we are interested in the case $c_1=1$, in this paper we will use $\beta = -1/2$ and $\alpha = -d/4$.  Notice that this algorithm will require the intrinsic dimension $d$ of the manifold $\mathcal{M}$, however we will determine this empirically as part of the kernel density estimation of the sampling density $q$.  


To determine the sampling density $q$ (which also serves as an estimate of the invariant measure $\peq$) we introduce the ad-hoc bandwidth function $\rho_0(x) = \left(\frac{1}{k_0-1}\sum_{j=2}^{k_0} ||x_i-x_{\textup{I}(i,j)}||^2 \right)^{1/2}$, where $\textup{I}(i,j)$ is the index of the $j$-th nearest neighbor of $x_i$ from the data set.  Following \cite{BH14VB} we used $k_0=8$ in all of our examples, and empirically the algorithm is not very sensitive to $k_0$.  With this bandwidth we define the following kernel,
\begin{align}\label{adhockernel} K_{\epsilon}(x,y) =  \exp\left(-\frac{||x-y||^2}{2\epsilon \rho_0(x)\rho_0(y)}\right), \end{align} 
which will be used only for the kernel density estimate of the sampling density $q$.  The kernel density estimate $q_{\epsilon}$ of the sampling density $q$ is given by the standard formula,
\begin{align}\label{KDE} q_{\epsilon}(x_i) &\equiv \frac{1}{N(2\pi\epsilon \rho_0(x_i)^2 )^{d/2}} \sum_{j=1}^N K_{\epsilon}(x_i,x_j)  \\
&\approx \int_{\mathcal{M}} \frac{K_{\epsilon}(x_i,y)}{(2\pi\epsilon \rho_0(x_i)^2 )^{d/2}}q(y)\, dV(y) = q(x_i) + \mathcal{O}(\epsilon), \nonumber \end{align}
where $dV$ is the volume form which $\mathcal{M}$ inherits from the ambient space and $q$ is the sampling measure relative to this volume form.  Note that applying \eqref{KDE} requires choosing the bandwidth $\epsilon$ and knowing the dimension $d$ of the manifold.

To determine the bandwidth $\epsilon$ and the dimension $d$, we apply the automatic tuning algorithm, originally developed in \cite{coifman2008TuningEpsilon} and refined in \cite{BH14VB}. The idea is that if $\epsilon$ is not well tuned, the kernel will become trivial; when $\epsilon$ is too small the kernel \eqref{adhockernel} is numerically zero when $x \neq y$, and when $\epsilon$ is too large the kernel is numerically one.  Forming the double sum $T(\epsilon) \equiv \frac{1}{N^2}\sum_{i,j=1}^N K_{\epsilon}(x_i,x_j)$, when $\epsilon$ is too small we find $T$ approaches $1/N$ and when $\epsilon$ is too large we find $T$ approaches $1$.  As shown in \cite{coifman2008TuningEpsilon,BH14VB} when $\epsilon$ is well tuned we have $T(\epsilon) \approx \frac{(4\pi\epsilon)^{d/2}}{\textup{vol}(\mathcal{M})}$ so that 
\begin{align}\label{Tofeps} \log T(\epsilon) \approx \frac{d}{2}\log(4\pi \epsilon) - \log(\textup{vol}(\mathcal{M})). \end{align} 
Since $T(\epsilon)$ is monotonically increasing in $\epsilon$, we also have $\log (T(\epsilon))$ monotonically increasing in $\log(\epsilon)$, and so the derivative $\frac{d \log(T(\epsilon))}{d\log(\epsilon)}$ has a unique maximum.  Intuitively this maximum corresponds to $\epsilon$ that gives the maximum `resolution' of the kernel $K_{\epsilon}$.  The approximation \eqref{Tofeps} suggests that the value of the maximum is $\max_{\epsilon} \frac{d \log(T(\epsilon))}{d\log(\epsilon)} = \frac{d}{2}$, and we will use this to determine the intrinsic dimension of the manifold $\mathcal{M}$.  Since the summation $T(\epsilon)$ is not very expensive to compute, we simply evaluate $T(\epsilon)$ for $\epsilon = 2^{l}$ where $l = -30,-29.9,...,9.9,10$ and then compute the empirical derivative,
\begin{align} \frac{d(\log S)}{d(\log \epsilon)} \approx \frac{\log(S(\epsilon+h))-\log(S(\epsilon))}{\log(\epsilon+h)-\log(\epsilon)}, \end{align}
and choose the value of $\epsilon$ which {\it maximizes} this derivative and set $d = 2\max_{\epsilon}\frac{d \log(T(\epsilon))}{d\log(\epsilon)}$.  

Now that we have the empirical estimate $q_{\epsilon}(x_i) = q(x_i) + \mathcal{O}(\epsilon)$, we can form the kernel \eqref{vbkernel}.  We reapply the same bandwidth selection to choose the global bandwidth $\epsilon$ in the kernel \eqref{vbkernel} and a new estimate of $d$, by constructing $T(\epsilon)$ as a double sum of the kernel \eqref{vbkernel} over the data set. Notice that these new values of $\epsilon$ and $d$ can be different from the previous values used in the kernel \eqref{KDE}, although empirically the new dimension $d$ is typically very similar.  

We can now evaluate the kernel \eqref{vbkernel} on all pairs from the data set and form the matrix $K^S_{\epsilon,i,j} = K^S_{\epsilon}(x_i,x_j)$ and we can compute the first normalization factor $q_{\epsilon}^S(x_i)= \sum_{j=1}^N K^S_{\epsilon}(x_i,x_j)/q_{\epsilon}(x_i)^{d\beta}$ as in Corollary \ref{corollary}.  We define a diagonal matrix $D_{i,i} = q_{\epsilon}^S(x_i)$, and the first normalization is to form the matrix $K^S_{\epsilon,\alpha} = D^{-\alpha} K_{\epsilon}^S D^{-\alpha}$.  We then compute the second normalization factor $q_{\epsilon,\alpha}^S(x_i) = \sum_{j=1}^N K^S_{\epsilon,\alpha,i,j}$ and form a diagonal matrix $D_{\alpha,i,i} = q_{\epsilon,\alpha}^S(x_i)$.  The second normalization is to form the matrix $\hat K_{\epsilon,\alpha}^S = D_{\alpha}^{-1}K_{\epsilon,\alpha}^S$.  We define the final normalization diagonal matrix $\hat{D}_{i,i} = 2 \epsilon q(x_i)^{2\beta}$, and by Corollary \ref{corollary},
\begin{align} L_{\epsilon,\alpha,\beta}^S = \hat D^{-1}(\hat K_{\epsilon,\alpha}^S - \textup{Id}) = \hat D^{-1}D_{\alpha}^{-1}D^{-\alpha} K_{\epsilon}^S D^{-\alpha} - \hat D^{-1} \nonumber \end{align}
approximates the desired operator $\hat{\mathcal{L}}$ when $\beta=-1/2$ and $\alpha =-d/4$.  To find the eigenvectors of $L_{\epsilon,\alpha,\beta}^S$, which approximate the eigenfunctions $\varphi_j$ of $\hat{\mathcal{L}}$, we note that setting $P = \hat D^{1/2}D_{\alpha}^{1/2}$ we can define a symmetric matrix,
\[ \hat{L} \equiv P L_{\epsilon,\alpha,\beta}^S P^{-1} = P^{-1}D^{-\alpha}K_{\epsilon}^S D^{-\alpha}P^{-1} - \hat D^{-1}. \]
Since $\hat{L}$ is a symmetric matrix, which is conjugate to the $L_{\epsilon,\alpha,\beta}^S$, we can compute the eigenvectors of $\hat{L} = \hat U\Lambda \hat U^{\top}$ efficiently and then the eigenvectors of $L_{\epsilon,\alpha,\beta}^S$ are given by the column vectors of $U = P^{-1}\hat U$.  

Note that the columns of $\hat U$ will be numerically orthogonal, so the columns of $U$ are orthogonal with respect to $P^2$ since $\textup{Id} = \hat U^{\top} \hat U = U^{\top}P^2 U$.  A careful calculation based on the asymptotic expansions in \cite{BH14VB} shows that $P_{ii}^2 = q(x_i)^{c_1-1} + \mathcal{O}(\epsilon)$ and in general the $q^{c_1}$ is the invariant measure of the gradient flow \eqref{gradSDE} so that $P^2$ represents the ratio between the invariant measure $e^{-c_1 U}$ of \eqref{gradSDE} and the sampling measure $q$.  However, in this case since $c_1=1$, we have $P = \textup{Id} + \mathcal{O}(\epsilon)$.  Thus, for the case $c_1=1$, we will take the eigenvectors $\varphi_j$ to be the column vectors of $\hat U$, since these eigenvectors are numerically orthogonal and are equal to the column vectors of $U$ up to order-$\epsilon$.  Notice that the orthogonality of these vectors $\frac{1}{N}\sum_{i=1}^N \varphi_l(x_i)\varphi_j(x_i) \approx \left<\varphi_l,\varphi_j\right>_{q}$, corresponds to the orthogonality of the eigenfunctions $\varphi_j$ with respect to the sampling measure (and since $c_1=1$ this is also the invariant measure of \eqref{gradSDE}).  Finally, in order to insure that the eigenvectors are orthonormal, we renormalize each column vector so that $\frac{1}{N}\sum_{i=1}^N \varphi_j(x_i)^2 = 1$.

\section{Appendix B: Estimating the Semi-group solutions for non-elliptic operators}\label{buildmodel}

Consider a dynamical system for a state vector $x$ on a manifold $\mathcal{M} \subset \mathbb{R}^n$ given by, 
\begin{align}\label{fullSDE} dx = a(x)\, dt + b(x)\, dW_t \end{align}
where $W_t$ is a standard Brownian process (generated by the Laplacian on $\mathcal{M}$), $a$ is a vector field (which is \emph{not necessarily} the gradient of a potential), and $b$ a diffusion tensor, all defined on $\mathcal{M}$.  Let $x_i = x(t_i)$ be a time series realization of \eqref{fullSDE} at discrete times $\{t_i\}_{i=1}^{N+1}$ with $\tau = t_{i+1}-t_i$.  As in the previous section, we represent a smooth function $f$ on the manifold by a vector $f_i = f(x_i)$.  Using the time ordering of the data set, we can define the shift map, $Sf(x_i) = f(x_{i+1})$.  Applying the It\^o formula to the process $y(t) = f(x(t))$ we have,
\begin{align}\label{generator}  
dy(s) &= \left( a\cdot \nabla f +  \frac{1}{2}\textup{Tr}(b^\top \textup{H}(f) b) \right)\, ds + \nabla f^\top b\, dW_s \nonumber\\ &\equiv \mathcal{L}f \, ds + \nabla f ^\top b\,dW_s, 
\end{align}
where $\textup{H}(f)$ denotes the Hessian and the functions and derivatives are evaluated at $x(s)$.  We first show that the expected value of the discrete time shift map is the semi-group solution $e^{\tau\mathcal{L}}$ associated to the generator $\mathcal{L}$ of the system \eqref{fullSDE}.  

For all smooth functions $f$ defined on the data set, the shift map yields a function $Sf$ which is defined on the first $N-1$ points of the data set.  Rewriting \eqref{generator} we have,
\begin{align}\label{shift1} Sf(x_i) &= f(x(t_{i+1})) = y(t_{i +1}) \\ 
&= f(x_i) + \int_{t_i}^{t_{i+1}} \mathcal{L}f \, ds + \int_{t_i}^{t_{i+1}} \nabla f^\top b\, dW_s, \nonumber \end{align}
 and taking the expectation conditional to the state $x_i$,
 \begin{align}\label{conditional1} \mathbb{E}_{x_i}[Sf(x_i)] = f(x_i) + \int_{t_i}^{t_{i+1}} \mathbb{E}_{x_i}[\mathcal{L}f(x_s)] \, ds. \end{align}
Recall that by the Feynman-Kac connection, the conditional expectation of the functional $y(t_{i+1}) = f(x_{i+1})$ is given by the semi-group solution $\mathbb{E}_{x_i}[y(t_{i+1})] = e^{\tau\mathcal{L}}f(x_i)$ and combining this with \eqref{conditional1} we find,
\begin{align}\label{conditional2} e^{\tau\mathcal{L}}f(x_i)  = f(x_i) + \int_{t_i}^{t_{i+1}} \mathbb{E}_{x_i}[\mathcal{L}f(x_s)]. \end{align}
Substituting \eqref{conditional2} into \eqref{shift1} we find,
\begin{align}\label{conditional3} Sf(x_i) &=  e^{\tau\mathcal{L}}f(x_i)  + \int_{t_i}^{t_{i+1}} \nabla f^\top b\, dW_s \nonumber\\ &\hspace{10pt}+ \int_{t_i}^{t_{i+1}} \mathcal{L}f - \mathbb{E}_{x_i}[\mathcal{L}f] \, ds. \end{align}
The formula \eqref{conditional3} shows that the expectation of the shift map $S$ is the semi-group solution $e^{\tau\mathcal{L}}$ as claimed.  This suggests that we can use the shift map to estimate the semi-group solution of the generator of \eqref{fullSDE}.  We next show that, by representing $S$ in an appropriate basis, we can minimize the error of this estimate. 

From \eqref{conditional3}, for any smooth function $g$, we can define the Monte-Carlo integral,
\begin{align}\label{mc} \left<g, Sf \right>_{\peq} &= \lim_{N\to\infty} \frac{1}{N-1}\sum_{i=1}^{N-1} g(x(t_i))Sf(x(t_i)) \nonumber  \\
&= \left<g,e^{\tau\mathcal{L}}f\right>_{\peq} + \left<g, \int_{t_i}^{t_{i+1}} \nabla f ^\top b \, dW_s \right>_{\peq} \nonumber \\ &\hspace{10pt}+ \left<g, \int_{t_i}^{t_{i+1}} Bf \, ds \right>_{\peq}, 
\end{align}
where we define $Bf = \mathcal{L}f - \mathbb{E}_{x_i}[\mathcal{L}f]$.  The Monte-Carlo integral implies that the inner products should be taken with respect to the sampling measure for the training data set, and we assume that the evolution of $x$ is ergodic so that the sampling measure is the invariant measure $\peq$ of the system \eqref{fullSDE}.  
Note that for smooth functions $f,g \in L^2(\mathcal{M},\peq)$, the final integral in \eqref{mc} will be order-$\tau$ since it is deterministic and the inner product with $g$ will be bounded.  Therefore, our goal is to choose $f, g$ from an orthonormal basis for $L^2(\mathcal{M},\peq)$ which minimizes the inner product with the stochastic integral, and thereby reduces the variance of our estimates of the coefficients.  

We first expand the norm of $\Omega(f) = \int_{t_i}^{t_{i+1}} \nabla f ^\top b \, dW_s$ by applying the It\^o isometry,
\begin{align} \left|\left|\Omega(f) \right|\right|_{\peq}^2 &= \int_{\mathcal{M}}\left(\int_{t_i}^{t_{i+1}} \nabla f ^\top b \, dW_s \right)^2 \peq(x_i) \, dV(x_i) \nonumber \\ &= \int_{\mathcal{M}}\int_{t_i}^{t_{i+1}} (\nabla f ^\top b)^2 \, ds \,  \peq(x_i) \, dV(x_i) \nonumber \\
&= \tau ||\nabla f^\top b ||_{\peq}^2 + \mathcal{O}(\tau^2).
\end{align}
In order to have a simple generic approach we will avoid estimating the diffusion tensor $b$ by assuming that the norm $||b(x)|| = \sup \frac{||v^\top b||}{||v||}$ is bounded above on the manifold by a constant $b_0$.  We can now apply the Cauchy-Schwartz inequality to find that,
\begin{align} |\left<g, \Omega(f) \right>_{\peq} | &\leq ||g||_{\peq}||\Omega(f)||_{\peq}\nonumber \\ &\leq \sqrt{\tau} b_0 ||g||_{\peq}||\nabla f||_{\peq} + \mathcal{O}(\tau). \end{align}
The optimal basis will be the one that minimizes the norm $||\nabla f||_{\peq} = \int_{\mathcal{M}} |\nabla f |^2 \peq dV(x).$
As shown in Appendix A, the norm $||\nabla f||_{\peq}$ is provably minimized by the eigenfunctions of the generator $\hat{\mathcal{L}}$ of the gradient flow system \eqref{gradSDE} and these eigenfunctions can be estimated by the diffusion kernel \eqref{vbkernel}.  Letting $\varphi_i$ be the eigenfunctions of $\hat{\mathcal{L}}$ with eigenvalues $\lambda_i$ so that $\hat{\mathcal{L}}\varphi_i = \lambda_i\varphi_i$, the set $\{\varphi_i\}$ is a basis for $L^2(\mathcal{M},\peq)$ and the norm we wish to minimize is given by the eigenvalues $||\nabla \varphi_i||_{\peq} = \lambda_i$.  

Replacing $f$ and $g$ in \eqref{mc} with these eigenfunctions, for a finite data set we define $\hat A_{lj} \equiv \frac{1}{N} \sum_{i=1}^N \varphi_j(x_i)S\varphi_l(x_i)$ and $A_{lj} \equiv \left<\varphi_j,e^{\tau\mathcal{L}}\varphi_l\right>_{\peq}$.  Since $\mathbb{E}[\varphi_j(x_i)S\varphi_l(x_i)]  = A_{lj}$, we have $\mathbb{E}[\hat A_{lj}] = A_{lj}$ which shows that $\hat A_{lj}$ is an unbiased estimate of $A_{lj}$.  Using the error bounds derived above, the variance is,
\begin{align} \mathbb{E}[(\hat A_{lj} - A_{lj})^2] 
&\leq b_0^2\lambda_l^2 \tau N^{-1} +\mathcal{O}(\tau^2 N^{-1}),  \nonumber \end{align}
assuming that $x_i$ are independent. Since $x_i$ form a time series, they are not independent and the convergence of the Monte-Carlo integral will be slower if the dependence is strong. In that case, one may need to subsample the time series which simultaneously requires a larger data set. Assuming independence, by the Chebyshev bound,
\begin{align} P( |\hat A_{lj} - A_{lj} | \geq \epsilon ) 
&\leq b_0^2\lambda_l^2\tau \epsilon^{-2}N^{-1} + \mathcal{O}(\tau^2 \epsilon^{-2}N^{-1}) \nonumber
\end{align}
and balancing these error terms requires $\lambda_l < b_0^{-1}\sqrt{\tau}$ and the errors are of order $\epsilon = \mathcal{O}(\tau N^{-1/2})$ in probability.

\end{document}